\documentclass[12pt, twoside,a4paper]{article}
\pdfoutput=1
\usepackage[cp1251]{inputenc}
\usepackage{ifthen}
\usepackage{amsmath}
\usepackage{amsfonts}
\usepackage{latexsym}
\usepackage{amsthm}
\usepackage{amssymb}
\newcommand{\CopyName}{ V.\ M.\ Zhuravlov}
\newcommand{\NAME}{ V.\ M.\ Zhuravlov}
\newcommand{\Year}{2024}
\newcommand{\rightheadtext}{Uniform terms and local elements}
     \newcounter{chapter}
     \newcounter{artpage}[chapter]
     \newcommand{\vs}{\vspace{.1in}}
     \newcommand{\vsk}{\vspace{.2in}}
\makeatletter
     \renewcommand{\@evenhead}{\footnotesize \ifthenelse{\value{artpage}=0}
     {\hfil}{\thepage\hfil \textsc {\leftmark} \hfil } }
     \renewcommand{\@oddhead}{\footnotesize\ifthenelse{\value{artpage}=0}
     {\hfil}{\hfil \textsc \rightmark \hfil \thepage} }
     \newcommand{\logo}{\baselineskip2pc \hbox to\hsize{\hfil\copyright\,\footnotesize
     \CopyName, \Year}}
     \renewcommand{\@oddfoot}{\ifthenelse{\value{artpage}=0}{\logo
     \refstepcounter{artpage}} {\hfil\refstepcounter{artpage}}}
     \renewcommand{\@evenfoot}{\ifthenelse{\value{artpage}=0}{\logo
     \refstepcounter{artpage}} {\hfil\refstepcounter{artpage}}}
     \renewcommand{\section}{\@startsection{section}{1}{0pt}{3.5ex plus
     1ex minus .2ex}{2.3ex plus 2.ex}{\large\hfil\textsc}}
     \vspace{3pt}\vs

\vs
\newcommand{\tit}{Uniform terms and local elements}
\date{2024}
\usepackage{hyperref}
\usepackage{graphicx}
\graphicspath{{Images/}}
\begin{document}
\hfill
\vspace{0.3in}
\markboth{{\NAME}}{{\rightheadtext}}\begin{center} \textsc {\CopyName} \end{center}\begin{center} \renewcommand{\baselinestretch}{1.3}\bf {\tit} \end{center}
\vspace{20pt plus 0.5pt} {\abstract{\noindent
The article explores function terms within uniform theories. It examines the uniformity of these theories through an algebraic lens. The paper compares the uniformity of terms and predicates within axiom schemas. It demonstrates the connection between such theories and both field theory and measure theory. The work proposes a novel interpretation of elements as automorphisms of the sets that include them. Additionally, it introduces a method for the localization of elements in relation to one another.\newline
\textit{Uniform Terms and Local Elements, 2024, msc: 03G30, 03H99\vspace{3pt}}\newline
}
}\vsk
\tableofcontents
\section{Introduction}
This article is a continuation of the paper: \href{https://arxiv.org/abs/2307.00069}{"Uniformity and Nonuniformity."} Specifically, we will explore the uniformity of the theory's terms and, consequently, the uniformity of operations and a new approach to localizing elements within an individual domain, or even constructing non-singleton theories. A term is defined as a predicate with a corresponding functional property, that is, with the property of uniqueness of its value, depending only on the values of its variables. Thus, we are interested in the algebraic properties of uniform and non-uniform sets.\par
Let us recall the definition of uniformity as given in the aforementioned article. Namely, a theory is \textbf{1}-uniform if the following axiom scheme holds true:
$$\exists x:A(x)\Longrightarrow\forall x:A(x)$$
— where \textbf{A(x)} is any unary logical formula without constants. And the theory is called \textbf{n}-uniform if the following scheme is true:
$$\exists(x_{1},...,x_{n}):\bigwedge_{j,k\in(1...n)}(x_{j}\neq x_{k})\wedge A(x_{1},...,x_{n}) \Longrightarrow \forall(x_{1},...,x_{n}):A!(x_{1},...,x_{n})$$
— where the quantifiers from a sequence imply a sequence of quantifiers (refer to the cited article above: \href{https://arxiv.org/pdf/2307.00069}{https://arxiv.org/pdf/2307.00069})\newline
\textbf{A!} denotes that there exists a permutation of variables for which the formula \textbf{A} holds true.\par
Initially, let's consider finitely uniform sets. These could be rationally or continually ordered lines or circles with a binary or ternary linear order. Fixing a certain constant within such a set divides it into two distinct parts (elements less than or greater than the constant; in cyclic order, things appear differently, but the differences are negligible), excluding the constant itself. Any element can serve as this distinguished constant. Moreover, any non-empty segment between a pair of distinct elements will, as we have already noted, be isomorphic to any other segment. Concurrently, there exists a group of automorphisms of the entire set \textbf{\textit{onto}} itself. It seems impossible to identify subgroups of movements and parallel translations within this group (which would be a normal subgroup). To make such an identification, it is necessary to introduce a metric on our ordered set—yet, we do not have any metric at present. It is only clear that parallel translations (shifts) will constitute an Abelian group, all transformations of which will either elevate or lower (in the case of transformations inverse to the elevating ones) the ordinal position of the elements of the entire set. However, this is not the only clarity. Shifts can be defined not only through a metric but also through the "external" behavior of an object within a certain class of other objects. There is no need to derive shifts from distance if it is possible to proceed conversely.
\section{Uniformity as a property of homomorphisms and automorphisms}
Shifts move all elements uniformly, without fixed points (or with all points being fixed). When the shift \textbf{f} is applied again to some element \textbf{c}, we get: $ff(c) \equiv f(f(c))$. In this case, the interval $[c,f(c))\equiv\{x|c\leq x<f(c)\}$, which is closed on the left and open on the right, will entirely shift into the interval $[f(c), ff(c))\equiv\{x|f(c)\leq x<ff(c)\}$. The reverse shift would set the corresponding left movement. These intervals have an empty intersection and a "complete" union, in the sense that $[c, f(c))\cup[f(c), ff(c))=[c, ff(c))$. A few words about the Abelian property of the shift group. There are a number of well-known, but very important properties behind the commutativity of operators. For any function (or operator), its graph is a set of pairs, which I prefer to call its construct (since a graph is, after all, a visual model of the structure created by the operator): $\{(x, f(x)) | x \in L\}$. Under the influence of the operator \textbf{g}, the \textbf{f}-construct will transform into $\{(g(x), gf(x)) | x \in L\}$. And it is quite obvious that the \textit{operators f and g commute if and only if one of them preserves the construct of the other, which is also equivalent to the operators not acting (leaving them fixed) on each other: $f = g^{-1}fg$ and $g = f^{-1}gf$.} This is why in quantum physics certain quantities have definite values simultaneously if and only if their operators are commutative.\par
Let's return to our discussion of shifts. We consider that a finitely uniform set \textbf{L} belongs to a certain algebraic universe (whether it be a universal set (universal class), or a category of ordered sets), "encompassing" a collection of homomorphisms sufficient for our discourse — this might sound a bit nebulous, but it is intuitively understandable. In the category of partially ordered sets and their homomorphisms, there exists a natural number object \textbf{N} with its own 'successor' function—the operator $\hat{s} \equiv (+1)$. This object enumerates the iterations of any arrows in the category, including the shifts of a rationally or continually ordered line; that is, all squares of the following diagram are commutative:\begin{center}
\includegraphics{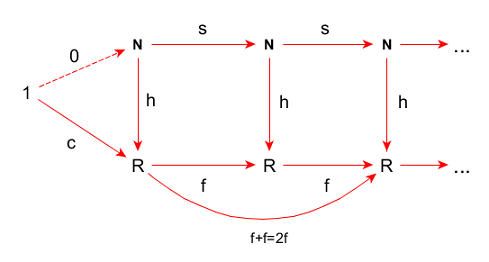}NR)\end{center}
The diagram operates 'to the left' as well, for the inverse shifts: \textbf{(-f)}. For every shift, there exists an interval that is transformed into the 'immediately following' interval. The NR diagram defines a measure on the set \textbf{R} — initially, all intervals are equipotent. The set of all shifts of the set \textbf{L} is an ordered set, isomorphic to \textbf{L} with an arbitrary constant point corresponding to the identity shift. Each shift congruently moves all intervals, thus preserving the constructs of all other shifts. The shifts of \textbf{L} are as uniform as \textbf{L} itself. For these reasons, the shifts are commutative. Furthermore, the diagram below is not just commutative but also Cartesian:\begin{center}
\includegraphics{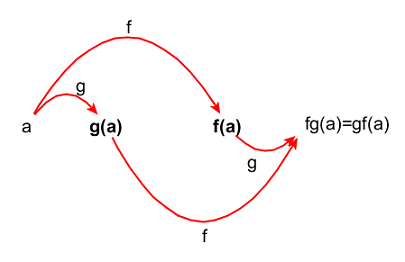}\end{center}
— the main thing here is the internal congruences of the segments:
$$\hat{g}[a, f(a)]=[g(a), gf(a)]$$ $$\hat{f}[a, g(a)]=[f(a), fg(a)]$$\par
In addition, there are the following universal algebraic properties of divisibility:
$$\forall x\forall s\exists !h:(x=s\cdot h)$$ $$\forall x\forall s\exists !h:(x=h\cdot s)$$\begin{center}
\includegraphics{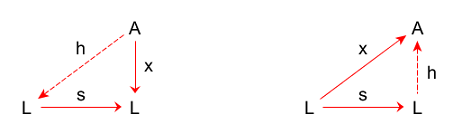}\end{center}
— where \textbf{x} is an arbitrary morphism, and \textbf{s} is an arbitrary shift of \textbf{L}. According to the second of these diagrams, \textit{every shift is an isomorphism, as the identity $1_{L}$ is uniquely decomposable into the product of an arbitrary shift and its inverse.}\par
And most importantly: instead of the categorical axiom \textbf{NR}, which asserts the existence of a natural number object, we adopt the statement \textbf{Z} about the existence of a (non-singleton) integer object:\begin{center}
\includegraphics{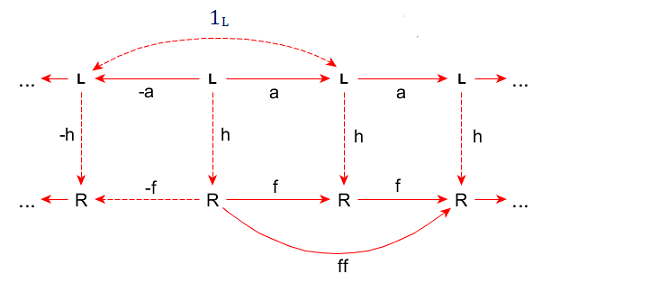} Z)\end{center}
— where all squares are Cartesian, the morphisms \textbf{h} and \textbf{(-h)} are unique; in addition, for every homomorphism \textbf{f} there is a homomorphism \textbf{(-f)} whose lifting along \textbf{(-h)} gives the pair \textbf{(-a, h)}.
From axiom \textbf{Z)}, in particular, it follows that squares are Cartesian:\begin{center}
\includegraphics{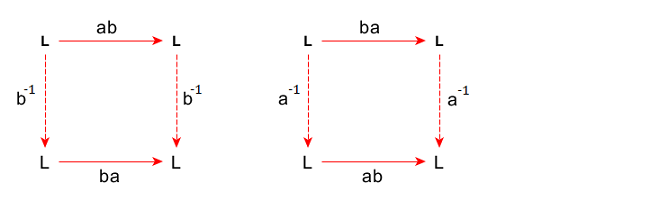}\end{center}
— from which, it seems to me, follows the Abelian property of the group of shifts.\par
The condition of a linear order being densely populated everywhere on \textbf{L} is equivalent to the condition of the morphism \textbf{a} being infinitely divisible (there are infinitely many such morphisms \textbf{a} that satisfy the \textbf{Z)}-property).\par
We have intentionally omitted the consideration of elements in the \textbf{Z}-axiom: for absolutely uniform \textit{non-standard sets}, the presence of their elements is not necessary. At the same time, a certain family of monomorphisms can always serve as elements.

\section{The group of shifts from the perspective of measure theory}
Let us now write down such a (possibly redundant) definition of a shift.\newline
\textbf{Definitions.} \textit{In the group of automorphisms of an everywhere dense (in the rational or continual sense), unbounded linearly ordered set \textbf{L} onto itself (i.e., a set that is finitely or continuously uniform), there is an Abelian subgroup \textbf{G} with the following properties:}
$$\forall f:(\forall x: f(x)<x)\bigvee\forall x:(f(x)=x)\bigvee(\forall x: f(x)>x)$$

$$\forall(f,c):c<f(c)\Longrightarrow[c,f(c))\bigcap[f(c),f^{2}(c))=\emptyset$$
$$\forall(f,c):c<f(c)\Longrightarrow[c,f(c))\bigcup[f(c),f^{2}(c))=[c,f^{2}(c))$$

$$\forall(f,c):c>f(c)\Longrightarrow[f(c),c)\bigcap[f^{2}(c),f(c))=\emptyset$$
$$\forall(f,c):c>f(c)\Longrightarrow[f(c),c)\bigcup[f^{2}(c),f(c))=[f^{2}(c),c)$$
\textit{— where: $ f^2 \equiv f\cdot f $; thus, it turns out that the statement about covering the entire \textbf{L} with segments specified by integer powers of any translation is equivalent to Archimedes' axiom. Let's also define:}
$$f\geq g \equiv \forall x:(f(x)\geq g(x))$$
\textit{— due to uniformity, in fact, one such \textbf{x} is enough.} This will be the group of translations of a finitely uniform set \textbf{L}.\par
If we now associate an arbitrary element $0 \in L$ with an identical translation, then each translation $f \in G$ can be associated with a single element $f(0) \in L$. Then it is easy to carry out the process of assigning coordinates, completely analogous to that used, for example, when constructing coordinates on the projective line. The mapping: $f \longmapsto f(0)$ is an order isomorphism with \textbf{L}. We also see that translations are exactly those automorphisms that preserve the additive measure on the set; and now we have shown how this measure is induced by successive executions of a translation (we must not forget about the automorphisms of the set of translations, which already define a multiplicative algebraic structure).\par
As a result, we get: axiom \textbf{Z)} defines a homomorphism of the set of integers into a finitely uniform set \textbf{L}. But \textbf{L} is dense everywhere. Therefore, the conjugate to this homomorphism is an isomorphism between \textbf{L} and the set of integers or real numbers. Which leads to specifying an additive measure function on \textbf{L}. And shifts (translations) are automorphisms that preserve the measure. Shifts define stationary disjoint coverings \textbf{L} (which transform into themselves when the shift is performed).\par
Thus, an everywhere dense linear order on \textbf{L}, unbounded from below and above, defines an ordered group (and even a field, if we take into account the automorphisms of the group) of its shifts, order isomorphic to \textbf{L}; it is enough just to identify each element with a shift that takes an arbitrarily chosen zero element to it. This is finite uniformity from an algebraic point of view.
\section{Field structure of finitely uniform sets}
\textbf{Theorem.} \textit{Terms of arity \textbf{2} and less, defined everywhere in a finitely uniform set with a binary predicate, define a field on this set. This field is unique, up to isomorphism.}\par
\textit{If a predicate on a uniform set is ternary, then the corresponding field is complemented by an element at infinity; thus, the coordinates of projective geometry will be determined on the set.}\par
\textit{The converse statements are also true.}\newline
\textbf{Proof.} Terms are functional symbols, i.e., \textbf{(n+1)}-predicates with the condition of uniqueness of one of the elements of the sequence.\par
But we are dealing with a finitely uniform set of elements \textbf{L}. The existence of a stationary element implies the stationarity of all elements:
$$\exists x:(x=\tau(x))\Longrightarrow\exists \forall x:(x=\tau(x))$$
— for some term $\tau$.
Consequently, there is a unit automorphism. Further, if an element has a preimage, then all elements have it:
$$\exists x\exists y:(x=\tau(y))\bigwedge(y\neq x)\Longrightarrow\forall x\exists y:(x=\tau(y))\bigwedge(y\neq x)$$\par
The same is true for images other than the elements that generate them:
$$\exists x\exists y:(y=\tau(x))\bigwedge(y\neq x)\Longrightarrow\forall x\exists y:(y=\tau(x))\bigwedge(y\neq x)$$\par
Consequently, unary terms are isomorphisms and constitute a group. The group is Abelian, due to the preservation of constructions (see above).\par
Now recall that finitely uniform sets are precisely everywhere dense, unbounded, linearly ordered sets. Since all elements of such a set are comparable, their uniformity with respect to unary terms implies:
$$\forall x\forall y:(x>y)\Longleftrightarrow(\tau(x)>\tau(y))$$
or:
$$\forall x\forall y:(x<y)\Longleftrightarrow(\tau(x)<\tau(y))$$
— unless the term is non-trivial. That is, non-trivial unary terms $\tau$ are automorphisms or anti-automorphisms of \textbf{L}. If any automorphism (or anti-automorphism) of set \textbf{L} elevates (or demotes) a certain element in the order hierarchy, then similarly, all elements will be affected:
$$ \exists x : (x \lessgtr \tau(x)) \Rightarrow \forall x : (x \lessgtr \tau(x)) $$
(do not forget that the requirement of non-constancy in the definition of uniformity concerns only \textbf{0}-ary terms, i.e., \textit{singleton} constants). Thus, among the unary terms, there is an Abelian group of (anti)automorphisms of a linearly ordered set.\par
If we now select an arbitrary element of the set $ 0 \in L $, identifying it with the identity element of the Abelian group, then each element \textbf{f} will correspond to a unique $ f(0) \in L $. Conversely, every element of \textbf{L} will correspond to a unique shift from the Abelian group of automorphisms of \textbf{L}.\par
As for binary terms, the product of elements of an Abelian group induces a corresponding binary operation on elements of \textbf{L}: $( (a \ast b) \equiv \widehat{a} (b) = \widehat{b} (a) )$.\par
Proceeding in a similar way with the automorphisms of the considered Abelian group, we naturally arrive at the structure of a field on \textbf{L}.\par
Moreover, the converse is also true: not only does a linear order determine the unique (up to isomorphism) algebraic structure of a continuous free field on \textbf{L}, but also, conversely, by declaring some elements of such a field positive, one can linearly order \textbf{L}.\par
In the case of a finitely uniform set with a ternary predicate, we have a cyclic ordering. The process of its "algebraization" is well described in textbooks on projective geometry, for example, by Efimov (see bibliography). The theorem has been proven.\par
\textbf{Note.} The theorem can also be proven by alternative methods, as previously mentioned.
\section{Equalization and localization}\par
The process of identifying elements of a uniformly finite set with the field of its automorphisms can be generalized.\par
It is important to note what is usually overlooked when considering isomorphic objects. Only category theory initially draws attention to the "truth up to isomorphism" of the overwhelming majority of our statements, and to the uniqueness "up to isomorphism" of most mathematical objects — mathematics usually does not go beyond the trivial acknowledgment of such concepts as isomorphism or equality, assuming that no interesting conclusions follow. However, in reality, isomorphism indeed represents an uncertainty, the inability to distinguish isomorphic objects by our means. And, on the other hand, isomorphism asserts the possibility of attributing to an object an essence, initially completely unexpected by us. The truth is that we cannot distinguish a uniformly finite set from its group of shifts. Furthermore, there is an infinite iteration of such ordered groups — directed both ways — as there is no first, original ordered set. And we are compelled to identify this set with its group of shifts — to equate these structures to one another. This is the simple idea of equating (or metamorphism), which we will repeatedly utilize further. I have already written (concerning the square root of \textbf{2}, or the imaginary unit) that in mathematics, there are no fictitious or real objects. Everything that we can conceive is real, and contradiction is the only perceptible obstacle to thought. And everything that we can conceive is real because it inevitably pertains to reality, regardless of whether we consider these fantasies to exist... Therefore, of all human endeavors, it is mathematics that is most closely linked with both our empirical existence and the world of our inventions and fantasies — which has remained a very mysterious and incomprehensible phenomenon.\par
But let us provide a verbal explanation of our constructions. If we demand operational distinguishability for a uniformly finite set \textbf{L}, then we must consider operators that transform each element into any other. The existence of inverse and identity operators is also evident. By applying operators sequentially, we reach their composition (analogous to arrows in category theory), independent of the order of bracket placement (only the sequence of execution matters). Uniformity also requires the same effect of the operator on all elements, leading to the invariance of their constructs regarding mutual influences — i.e., the set of operators forms an Abelian group. This group is free and infinite, just as \textbf{L} itself — uniform and infinite. Due to the identical action of the order-constructing operators on \textbf{L}'s elements, for each pair of elements, there is a unique operator that transforms the lesser into the greater. Therefore, by choosing an arbitrary $0 \in L$, all other elements correspond uniquely to the operators that generate them from \textbf{0}. This is the idea of equating — identifying an object with the operator (or set of operators) that generates it. We are discussing a certain essential principle, the cognitive value of which lies in significantly simplifying the overall picture.\par
In any case, by fixing an arbitrary point as \textbf{0}, we enable addition in a uniformly finite set. At the same time, each element simultaneously acts as an addition operator, a certain iso-arrow of \textbf{L} onto itself; the group of shifts is Abelian, hence the operator functions equally from both the right and the left. Now it suffices to scale \textbf{L} somehow — for instance, by identifying an arbitrary non-zero shift with \textbf{1} — to transform the entire line into rational (real) numbers. All details of this process can be found, for example, in Efimov's book (see bibliography) — there, only the orderliness of the projective line is utilized to construct coordinate numbers on it. Besides addition, it is crucial for us to define all other operations on the line.\par
We have demonstrated that on a uniformly finite set, one can define an addition operation localized relative to an arbitrary choice of zero; the arbitrariness here is that additions defined by different zeros are easily and unambiguously transformed among themselves. Therefore, by equating the original binary ordered set to an ordered Abelian group, we single out element \textbf{0}, but any other element could serve as such. We have already spoken about this situation, calling it existentially (in)uniformity,— the particle "not" in parentheses is used intentionally — due to the ambiguous semantics of such cases. \textbf{0} exists, but we do not know where it is, and thus we can place it anywhere... Simultaneously, our model set \textbf{L} already becomes a scale, the position within which carries much more information than before. And each shift already defines a measure on \textbf{L}. Now it is important to note that the free Abelian group we have obtained also has its automorphisms.\par
The original Abelian group is free (being a structure induced by uniform finiteness), ordered (the order is consistent with the positive elements — those greater than zero), and does not have a finite set of generators (due to the everywhere density of the order). And this group can now be equated to a field by equating its non-zero elements to the automorphisms of the group itself. And again, for this we need to choose an arbitrary positive element in \textbf{L} as \textbf{1}. The procedure for such a transformation of a group into a ring (and in this case, into a field) is known. The very definition of a group automorphism leads to the distributivity of the multiplicative law with respect to the additive law, etc. So we come to the conclusion that the model of a uniformly finite binary algebra is the field of rational or real numbers with their natural ordering and arbitrarily chosen \textbf{0} and \textbf{1}. It is interesting that the order anti-isomorphism $\times(-1)$ appears here; which “algebraically” should not confuse us — conceptually, these are simply two states of order, both equally equivalent to uniformity. It is also interesting, as we did in the case of shifts (summation transformations), to provide an ordinal description of the multiplicative operations:
$$\widehat{a}(x)\equiv a\ast x$$
$$\forall(a,x,y):\widehat{a}(x+y)=\widehat{a}(x)+\widehat{a}(y)$$
$$\widehat{a}[0, 1)=[0, a)$$
— with all the ensuing order consequences; the operator of multiplication by an element is also an order isomorphism or anti-isomorphism, but not a shift, rather a stretch (or compression). And here, too, there is an element corresponding to the "fixed" operator $ (\hat{1}\times ) $, arbitrarily fixed in \textbf{L}. Thus, a uniformly finite set is equated to the field of rational or real numbers. The idea here is that we cannot distinguish a shift of a uniformity set from a "stretched" shift, i.e., multiplied by a certain value.\par
So, having fixed a pair of constant elements \textbf{(0,1)} in a uniformity set, we can fix all other elements, assigning them specific coordinates. In other words, without an absolute distinction between elements, we can always distinguish and localize them relative to an arbitrary pair of elements. One can also arrive at integers or natural numbers by fixing an arbitrary shift-monomorphism in the real numbers as an atomic element.\par
As we have already shown earlier, absolutely non-uniformity sets coincide with completely ordered sets — those in which all elements are distinguishable. But completely ordered sets are defined up to isomorphisms. For any completely ordered set, each of its right segments is isomorphic to it; moreover, such a set itself is a certain right segment of a set isomorphic to it. By combining this entire family of sets and supplementing the result with multiplicative automorphisms of extensions (and compressions), we immerse the non-uniformity set into an absolutely non-uniformity one, which also contains other non-uniformity sets. Apparently, there is some conjugacy here.\par
And we can finally answer the question posed at the beginning of the cited article \href{https://arxiv.org/pdf/2307.00069}{"Uniformity and Nonuniformity"}: why are numbers needed? It turns out that a maximally uniformity series of objects, potentially capable of distinguishing something, of optimally fixing some information, and of operating with this information, is equivalent to a binary, everywhere dense linear order. The latter, in turn, is equivalent to a numerical field (with additive and multiplicative operations). This is the cognitive meaning of such constructions. This is their objective meaning — it is with their help that all precise science is constructed.
\section{Uniformity and connectivity}
Finite uniformity is equivalent to an everywhere dense linear ordering in the binary case, or an analogous cyclic ordering in the ternary case. When we discussed equating, we asserted that an element of a set is indistinguishable from certain properties that define it (for example, from a shift that translates another element into this one). However, predicates that define properties determine subsets of an individual domain; these subsets, in turn, can also relate to orderliness if we are discussing infinite supremums. It is known that there is the following Galois correspondence (or order self-adjointness of \textbf{L}, if you prefer):\par
$X^{>}\equiv\{x|x > X\}$, $X^{<}\equiv\{x|x < X\}$ (then: $X^{> < >}=X^{>}$ and: $X^{< > <}=X^{<}$).\\par
For any bounded set \textbf{X}, pairs of sets: $(X^{<}, X^{(< >)})$ or in the dual version: $ (X^{>}, X^{(> <)})$, serve as Dedekind cuts on \textbf{L} (these pairs can serve as the definitions of cuts). For such pairs in the set of rational numbers, there are two types of cuts:
$$\exists! x:(X^{<}\leq x)\wedge(x\leq X^{< >})$$ or:$$\neg\exists x:(X^{<}\leq x)\wedge(x\leq X^{< >})$$,
moreover, in the field of rational numbers, such \textbf{x} exist for some bounded sets \textbf{X} (this is true for cuts defined by elements of \textbf{L}). Therefore, if we require \textbf{L} to be even more uniform than the finitely uniform field of rational numbers, we must require the existence of elements for all cuts. That is, we come to the field of real numbers, where not only elements can be uniform, but also some unary predicates (or bounded subsets of \textbf{L}).\par
But then the question arises: what are these subsets? The answer lies in the ideas of equating outlined above. These are the subsets that define certain elements — existing or non-existing in \textbf{L}. As a result, we find that a linear (or cyclic) order is infinitely uniform if and only if \textbf{L} is a topologically connected set (all cuts have singleton closures). Otherwise, \textbf{L} would be the sum of a certain number of disjoint open sets.
\section{References}
\noindent [1] Joseph R. Shoenfield. Oxford.Addison-Wesley Publishing Company,\par
\textit{Mathematical Logic}, (1967).\newline
 [2] Haskell B. Curry,  McGraw-Hill Book Company, INC.\par
\textit{Foundations of Manthematical logic}, (1984).\newline
 [3] Helena Rasiowa and Roman Sikorski,  Panstwowe Wydawnlctwo Naukowe Warszawa,\par
 \textit{The Mathemayics of Metamathematics}, (1963).\newline
 [4] Alonzo Church, Princeton University Press,\par
 \textit{Introduction to Mathematical Logic}, (1956).\newline
 [5] Efimov N.V. Moscow. Science,\par
 \textit{Higher geometry}, (1971)\newline
 [6] P.M. Kon, D. Reidel Publishing Company,\par
\textit{Universal Algebra}, (1961).\newline
 [7] R. Goldblatt, North-Holland Publishing Company Amsterdam New York Oxford\par
 \textit{TOPOI. The categorial analysis of logic}, (1979)\newline
 [8] D. Hilbert und P. Bernays. Springer-Verlag Berlin — Heidelberg — New York,\par
 \textit{Grundlagen der Mathematik. I}, (1968).

\end{document}